\documentclass{article}

\usepackage{amsmath}

\usepackage{amsfonts}

\newtheorem{theorem}{Theorem}

\newtheorem{corollary}[theorem]{Corollary}

\newtheorem{definition}[theorem]{Definition}

\newtheorem{lemma}[theorem]{Lemma}

\newtheorem{proposition}[theorem]{Proposition}

\newenvironment{proof}[1][Proof]{\noindent\textbf{#1.} }{\ \rule{0.5em}{0.5em}}

\begin{document}

\begin{center}
{\Large \textbf{Some Characterizations in Kenmotsu manifolds with a new connection}}\bigskip \bigskip

Mohammed Ali Alghamdi and O\u{g}uzhan Bahad\i r\bigskip \bigskip
\end{center}

\noindent \textbf{Classification of Mathematics Subject:} 53C15, 53C25, 53C40.
\


\medskip \noindent \textbf{Abstract.}  The present study initially identified the generalized symmetric connections $(\alpha,\beta)$ typed, which can be regarded as more generalised forms of quarter and semi-symmetric connections.  The quarter and semi-symmetric connections are obtained respectively particularly when $(\alpha,\beta)=(1,0)$ and $(\alpha,\beta)=(0,1)$ are taken into consideration. Taking that into account, a new generalized symmetric metric connection was attained upon Kenmotsu manifolds. In compliance with the new connection, some results were provided through calculation of tensors belonging to Kenmotsu manifold involving curvature tensor, ricci tensor, projective curvature tensor and coincircular curvature tensor. Ultimately, an original example was presented.
 \medskip

\section{Introduction}

A particular metric connection with a torsion different from zero was introduced by Hayden upon a Riemannian manifold \cite{Hayden}. The quarter symmetric connections, being more generalized form of semi-symmetric connections,were suggested by Golab upon a differentiable manifold \cite{4}. As for the present study, the definition below is presented by taking it a step further and generalising the quarter symmetric connections as well.

A linear connection is suggested to be generalized symmetric connection on condition that the torsion tensor connection is presented in the form as follows:
\begin{small}
\begin{equation}
T(U,V)=\alpha \{u(V)U-u(U)V\}+\beta \{u(V)\varphi U-u(U)\varphi Y\},   \label{Int-2}
\end{equation}
\end{small}
for all vector fields $U$, $V$ upon a manifold in which $\alpha$ and $\beta$ are smooth functions. $\varphi$ can be viewed as a tensor of type $(1,1)$ and $u$ is regarded as a $1-$ form connected with the vector field which has a non-vanishing smooth non-null unit. Furthermore, the connection mentioned is suggested to be a generalized metric one when a Riemennian metric g in M is available as $\overline{D}g=0$; or else, it is non-metric.

In the equation (\ref{Int-2}), if $\alpha=0$ ($\beta=0$), then the generalized symmetric connection is called $\beta-$ quarter-symmetric connection ( $\alpha-$ semi-symmetric connection), respectively. Additionally, the connection  is decreased to a semi-symmetric one and a quarter-symmetric one when $(\alpha,\beta)=(1,0)$ and $(\alpha,\beta)=(0,1)$ are selected. Thus, it can be suggested that generalizing semi-symmetric and quarter symmetric connections paves the way for a generalized symmetric metric connection. Those two connections are of great significance both for the study of geometry and applications in physics. For instance Pahan, Pal and Bhattacharyya  studied generalized Robertson--Walker space-time with respect to quarter-symmetric connection \cite{pah}. Furthermore, many authors investigated the geometrical and physical aspect of different spaces \cite{avik}, \cite{jin}, \cite{lee}, \cite{leee}, \cite{pahan}, \cite{qu}, \cite{zhang}.
The connections including semi-symmetric and semi-symmetric non-metric ones were investigated in a Kenmotsu manifold respectively by the authors \cite{Pra}, \cite{sukla}, \cite{si} \cite{ilk1}, \cite{ilk2}, \cite{zhao} and \cite{zulekha}.

In the present paper, we defined new connection for Kenmotsu manifold, generalized symmetric metric connection. This connection is the generalized form of semi-symmetric metric connection and quarter-symmetric metric connection. The preliminaries are presented in Section 2 and Section 3 illustrates generalized symmetric connection for a Kenmotsu manifold. As for Section 4, The curvature and Ricci tensors of a Kenmotsu manifold and its scalar curvature are calculated in relation to generalized metric connection. Besides, it is found that first Bianchi identity is provided and Ricci tensor is symmetric with respect to generalized metric connection of type $(-1,\beta)$ and $(\alpha,0)$. Section 5, indicates that it is a generalized $\eta -$ Einstein manifold with reference generalized symmetric metric connection when a Kenmotsu manifold looks $\phi-$ projectively  flat concerning generalized metric connection. Via Section 6, on the other hand, an expression is presented in the interest of concircular curvature tensor with respect to generalized metric connection. Finally, Section 7 provides an example verifying some results of Section 4 and Section 5.

\section{Preliminaries}
In order to call a differentiable $M$ manifold of dimension $n=2m+1$ as practically contact metric one, a $(1,1)$ tensor field, a contravariant vector field,  a $1-$ form and a Riemannian metric $g$ should be admitted, which satisfy
\begin{eqnarray}
\psi \xi&=&0,\\
\eta(\psi U)&=&0\\
\eta(\xi)&=&1, \label{2.1}\\
\psi^{2}(U)&=&-U+\eta(U)\xi, \label{2.2}\\
g(\psi U,\psi V)&=&g(U,V)-\eta(U)\eta(V), \label{2.3}\\
g(U,\xi)&=&\eta(U), \label{2.4}
\end{eqnarray}
for all the vector fields $U$, $V$ on $M$. When $g(U,\psi V)=\Phi(U,V)$ is written, the tensor field $\phi$ is viewed as an anti-symmetric $(0,2)$ tensor field \cite{Bla}. When a practically contact metric manifold  performs
\begin{eqnarray}
(\nabla_{U}\psi)V&=&g(\psi U,V)\xi-\eta(V)\psi U, \label{2.9}\\
\nabla_{U}\xi&=&U-\eta(U)\xi,
\end{eqnarray}
M is regarded as a Kenmotsu manifold in which the Levi-Civita connection of g is present \cite{Ken}.

The following relations can be observed in Kenmotsu manifolds \cite{Ken}:
\begin{eqnarray}
(\nabla_{U}\eta)V&=&g(\psi U,\psi V)\\
g(R(U,V)W,\xi)&=&\eta(R(U,V)W)=g(U,W)\eta(V)-g(V,W)\eta(U), \label{2.11}\\
R(\xi,U)V&=&\eta(V)U-g(U,V)\xi, \label{2.12}\\
R(U,V)\xi&=&\eta(U)V-\eta(V)U, \label{2.13}\\
R(\xi,U)\xi&=&U-\eta(U)\xi, \label{2.14}\\
S(U,\xi)&=&-(n-1)\eta(U), \label{2.15}\\
S(\psi U,\psi V)&=&S(U,V)+(n-1)\eta(U)\eta(V) \label{2.16}
\end{eqnarray}
for all vector fields  $U$, $V$ and $W$, in which $R$ and $S$ can be viewed as the curvature and Ricci tensors belonging to $M$. A Kenmotsu manifold $M$ is found to be generalized $\eta$ Einstein when the Ricci sensor $S$ of it is in the form presented as follows;
\begin{eqnarray}
S(U,V)=ag(U,V)+b\eta(U)\eta(V)+cg(\psi U,V),
\end{eqnarray}
for every $U,V\in\Gamma(TM)$, in which $a$, $b$ and $c$ are viewed as scalar functions. In such a way that $b\neq0$ and $c\neq0$ in the event that $c=0$, M is regarded as $\eta-$ Einstein manifold.

\section{Generalized Symmetric Metric Connection in a Kenmotsu manifold}

When we view $\overline{D}$ as a linear connection and $D$ as a Levi-Civita connection of practically cantact metric manifold M in such a way that
\begin{eqnarray}
\overline{D}_{U}V=D_{U}V+H(U,V),\label{2.19}
\end{eqnarray}
for all vector field $X$ and $Y$. The following is obtained so that $\overline{D}$ is a generalized symmetric connection of $D$, in which $H$ is viewed as a tensor of type $(1,2)$;
\begin{eqnarray}
H(U,V)=\frac{1}{2}[T(U,V)+T^{'}(U,V)+T^{'}(V,X)],\label{2.20}
\end{eqnarray}
where $T$ is viewed as the torsion tensor of $\overline{D}$ and
\begin{eqnarray}
g(T^{'}(U,V),W)=g(T(W,U),V).\label{2.21}
\end{eqnarray}
Thanks to (\ref{Int-2}) and (\ref{2.21}), we obtain the following;
\begin{eqnarray}
T^{'}(U,V)=\alpha \{\eta (U)V-g(U,V)\xi\}+\beta\{-\eta (U)\phi V-g(\psi U,V)\xi\}.\label{2.22}
\end{eqnarray}
Furthermore, through use of (\ref{Int-2}), (\ref{2.20}) and (\ref{2.22}), we get
\begin{eqnarray}
H(U,V)=\alpha \{\eta (V)U-g(U,V)\xi\}+\beta\{-\eta (U)\psi V\}.,\label{2.23}
\end{eqnarray}
\begin{corollary}
For a Kenmotsu manifold, generalized symmetric metric connection $\overline{\nabla}$ is provided by
\begin{eqnarray}
\overline{D}_{U}V=D_{U}V+\alpha \{\eta (V)U-g(U,V)\xi\}-\beta\eta (U)\psi V.\label{konnek}
\end{eqnarray}
\end{corollary}

If  $(\alpha,\beta)=(1,0)$ and $(\alpha,\beta)=(0,1)$ are chosen, the generalized symmetric metric connection is diminished to a semi-symmetric metric and a quarter symmetric metric one as presented in the following;
\begin{eqnarray}
\overline{D}_{U}V=D_{U}V+\eta (V)U-g(U,V)\xi,\label{konnek1}
\end{eqnarray}
\begin{eqnarray}
\overline{D}_{U}V=D_{U}V-\eta (U)\psi V.\label{konnek2}
\end{eqnarray}
By means of (\ref{konnek}), we have the following proposition;
\begin{proposition} \label{propo}
The following relations are obtained when M is a Kenmotsu manifold with generalized metric connection:
\begin{eqnarray}
(\overline{D}_{U}\psi)V&=&(\alpha+1)\{g(\psi U,V)\xi-\eta (V)\psi U\},\\
\overline{D}_{U}\xi&=&(\alpha+1)\{U-\eta (U)\xi\},\\
(\overline{D}_{U}\eta)V&=&(\alpha+1)\{g(U,V)-\eta (V)\eta(U)\},
\end{eqnarray}
for every $U,V\in\Gamma(TM)$.
\end{proposition}

\section{Curvature Tensor}

Considering that $M$ is a $n-$ dimensional Kenmotsu manifold, the following can define the curvature tensor $%
\overline{R}$ of the generalized metric connection $\overline{D}$  on $M$.
\begin{eqnarray}
\overline{R}(U,V)W={\overline{D}}_{U}{\overline{D}}_{V}W-{\overline{D}}_{V}{\overline{D}}_{U}W-{\overline{D}}_{[U,V]}W, \label{3.1}
\end{eqnarray}
When proposition \ref{propo} is used, through (\ref{konnek}) and (\ref{3.1}), we obtain
\begin{eqnarray}
\overline{R}(U,V)W&=&R(U,V)W+\{(-\alpha^{2}-2\alpha)g(V,W)+(\alpha^{2}+a)\eta(V)\eta(W)\}U\ \nonumber\\
&&+\{(\alpha^{2}+2\alpha)g(U,W)+(-\alpha^{2}-\alpha)\eta(U)\eta(W)\}V \ \nonumber\\
&&+\{(\alpha^{2}+\alpha)[g(V,W)\eta(U)-g(U,W)\eta(V)]+(\beta+\alpha\beta)[g(U,\psi W)\eta(V)-g(V,\psi W)\eta(U)]\}\xi\ \nonumber\\
&&+(\beta+\alpha\beta)\eta(V)\eta(W)\psi U-(\beta+\alpha\beta)\eta(U)\eta(W)\psi V \label{eg}
\label{3.2}
\end{eqnarray}
where
\begin{eqnarray}
R(U,V)W={D}_{U}{\nabla}_{V}W-{\nabla}_{V}{\nabla}_{U}W-{\nabla}_{[U,V]}W, \label{3.111}
\end{eqnarray}
is viewed as the curvature tensor regarding Levi-connection $D$. When (\ref{3.2}) and the first Bianchi identity are paid attention, we have

\begin{eqnarray}
\overline{R}(U,V)W+\overline{R}(V,W)U+\overline{R}(W,U)V=2(\beta+\alpha\beta)\{\eta(U)g(\psi V,W)+\eta(V)g(U,\psi W)+\eta(W)g(V,\psi U)\}.
\end{eqnarray}
Thus, the following proposition is obtained
\begin{proposition}
Let us consider that $M$ is an n-dimensional Kenmotsu manifold together with generalized symmetric metric of type$(\alpha,\beta)$.
When  $(\alpha,\beta)=(-1,\beta)$ or $(\alpha,\beta)=(\alpha,0)$, it provides the first Bianchi identity of the generalized connection $\overline{D}$ on $M$.
\end{proposition}
The following proposition is presented using (\ref{2.11}), (\ref{2.12}), (\ref{2.13}), (\ref{2.14}) and (\ref{3.2}).
\begin{proposition}
When $M$ is an $n-$ dimensional Kenmotsu manifold with generalized symmetric metric of type $(\alpha,\beta)$, we have the following equations:
\begin{eqnarray}
\overline{R}(U,V)\xi&=&(\alpha+1)\{\eta(U)V-\eta(V)U+\beta[\eta(V)\psi U-\eta(U)\psi V]\} , \label{3.16}\\
\overline{R}(\xi,U)V&=&(\alpha+1)\{\eta(V)U-g(U,V)\xi+\beta[\eta(V)\psi U-g(U,\psi V)\xi]\} , \label{3.17}\\
\overline{R}(\xi,V)\xi&=&(\alpha+1)\{V-\eta(V)\xi-\beta \psi V\}, \label{3.18}\\
\eta(\overline{R}(U,V)W)&=&(\alpha+1)\{\eta(V)g(U,W)-\eta(U)g(V,W)+\beta[\eta(V)g(U,\psi W)-\eta(U)g(V,\psi W)]\}\ \nonumber\\
 \label{3.14}
\end{eqnarray}
for every $U,V,W\in\Gamma(TM)$.
\end{proposition}

In the following, the Ricci tensor $\overline{S}$ and the scalar curvature $\overline{r}$ of a Kenmotsu manifold is presented with generalized symmetric metric connection $\overline{D }$
$$\overline{S}(U,V)=\sum_{i=1}^{n}g(\overline{R}(\nu_{i},U)V,\nu_{i}),$$
$$\overline{r}=\sum_{i=1}^{n}\overline{S}(\nu_{i},\nu_{i}),$$
in which $U,V\in\Gamma(TM)$, $\{\nu_{1},\nu_{2},...,\nu_{n}\}$ is viewed as orthonormal frame. Then by using (\ref{2.3}) and (\ref{3.2}) we obtain

\begin{eqnarray}
\overline{S}(V,W)&=&S(V,W)+\{(2-n)\alpha^{2}+(3-2n)\alpha\}g(V,W)+(n-2)(\alpha^{2}+\alpha)\eta(V)\eta(W)\ \nonumber \\
&&-(\beta+\alpha\beta)g(V,\phi W), \label{Ricci}
\end{eqnarray}
in which $S$  is viewed as the Ricci tensor with regards to Levi-Civita connection.
Due to the fact that Ricci tensor of the Levi-connection is symmetric, (\ref{Ricci}) provides us the following
\begin{eqnarray}
\overline{S}(V,W)-\overline{S}(W,V)=-2(\beta+\alpha\beta)g(V,\psi W). \label{Riccim}
\end{eqnarray}
Thus, the theorem below is obtained.
\begin{theorem}
 Consider that $M$ is an $n-$ dimensional Kenmotsu manifold. The Ricci tensor $\overline{S}$ of generalized symmetric metric connection $\overline{D}$  is symmetric if and only  $(\alpha,\beta)=(-1,\beta)$ or $(\alpha,\beta)=(\alpha,0)$.
\end{theorem}
The following theorem is enabled through use of (\ref{Ricci})
\begin{theorem}
When $M$  is an $n-$ dimensional Kenmotsu manifold with generalized symmetric metric connection of type $(\alpha,\beta)$, The scalar curvature concerning the connection is as follows:
\begin{eqnarray}
\overline{r}=r+(n-2)(1-n)\alpha^{2}-2(n-1)^{2}\alpha, \label{3.122}
\end{eqnarray}
in which $r$ is viewed as scalar curvature of Levi-Civita connection.
\end{theorem}
\begin{lemma}
As in the case of Theorem 6, let us suppose that $M$ is an $n-$ dimensional Kenmotsu manifold with the connection, then we have
\begin{eqnarray}
\overline{S}(V,\xi)&=&(1-n)(\alpha+1)\eta(V), \label{3.12}\\
\overline{S}(\phi V,\phi W)&=&\overline{S}(V,W)+(n-1)(1+\alpha), \label{3.19}
\end{eqnarray}
for every $V,W\in\Gamma(TM)$.\\
\begin{proof}
Thanks to (\ref{2.1}), (\ref{2.15}) and (\ref{Ricci}), we have (\ref{3.12}). Anymore, the use of (\ref{2.3}), (\ref{2.16}) and (\ref{Ricci}) provides us (\ref{3.19}).
\end{proof}
\end{lemma}

\begin{theorem}
Suppose that $M$ is an $n-$ dimeansional Knemotsu manifold with the connection as in Lemma 7. $M$ is viewed as an $\eta-$ Einstein manifold regarding generalized symmetric metric connection $(for\;\beta \neq1\;and\; \beta\neq -1 )$ on condition that M is Ricci semi-symmetric regarding the connection $\overline{D}$.

\begin{proof}
Let $\overline{R}(U,V)\overline{S}=0$ be on $M$ for any $U,V,W,Z\in \Gamma(TM)$. Then the following is obtained;
\begin{eqnarray}
\overline{S}(\overline{R}(U,V)W,Z)+\overline{S}(W,\overline{R}(U,V)Z)=0.\label{et}
\end{eqnarray}
If we choose $W=\xi$ and $U=\xi$ in (\ref{et}), we get
\begin{eqnarray}
\overline{S}(\overline{R}(\xi,V)\xi,Z)+\overline{S}(\xi,\overline{R}(\xi,U)Z)=0 \label{etk}.
\end{eqnarray}
Using (\ref{3.17}), (\ref{3.18}) and (\ref{3.12}) in (\ref{etk}), we obtain
\begin{eqnarray}
\overline{S}(V,Z)-\beta\overline{S}(\phi V,Z)=(1-n)\{g(V,Z)+\beta g(V,\phi Z)\} \label{etkk}.
\end{eqnarray}
If one substitutes $V=\phi V$ in the equation (\ref{etkk}) and using (\ref{3.12}), we get
\begin{eqnarray}
\overline{S}(\phi V,Z)+\beta\overline{S}(V,Z)=(1-n)\{g(\phi V,Z)+\beta g(V,Z)+(\alpha-\beta+1)\eta(V)\eta(Z)\} \label{etkkk}.
\end{eqnarray}
From the (\ref{etkk}) and (\ref{etkkk}), we obtain
\begin{eqnarray}
\overline{S}(V,Z)=\frac{1-n}{1-\beta^{2}}\{((1+\beta^{2})g(V,Z)+(\alpha-\beta+1)\eta(V)\eta(Z)\} \label{etkkl}.
\end{eqnarray}
The equation suggests that $M$is an n Einstein manifold in regard to generalized metric connection.
\end{proof}
\end{theorem}
 \begin{corollary}
 Let us suppose that $M$ is an $n-$ dimensional Kenmotsu manifold with generalized symmetric metric connection of type $(\alpha,\beta)$.
  Considering M is Ricci semi-symmetric in regard to the connection, the statements illustrated below are obtained:\\
 \textbf{(i)} $M$  is reduced to an Einstein manifold with reference to generalized symmetric metric connection of type $(\beta-1,\beta)$.\\
 \textbf{(ii)} There is no Einstein manifold with the connection of type $(\alpha,1)$ and $(\alpha,-1)$.
 \end{corollary}

\section{Projective Curvature Tensor}

Let us consider that $M$ is an $n-$ dimensional Kenmotsu manifold with generalized connection of type $(\alpha,\beta)$. The following definition reveals the projective curvature tensor $P$ of type $(1, 3)$ of $M$ in regard to connection $\overline{\nabla}$:
\begin{eqnarray}
\overline{P}(U,V)W=\overline{R}(U,V)W-\frac{1}{n-1}\{\overline{S}(V,W)U-\overline{S}(U,W)V\}. \label{4.1}
\end{eqnarray}

\begin{definition}
Considering that $M$ is an $n-$ dimensional Kenmotsu manifold, $M$ can be suggested to be $\xi-$ projectively flat in regard to generalized symmetric metric connection $\overline{D}$ provided that $\overline{P}(U,V)\xi=0$ on $M$.
\end{definition}
The use of (\ref{3.16}), (\ref{3.12}) and (\ref{4.1}) enables us the following:
$$\overline{P}(U,V)\xi=(\alpha+1)\beta\{\eta(V)\psi U-\eta(U)\psi V\}.$$

Hence, the theorem presented below is obtained
\begin{theorem}
Let us suppose that $M$ is an $n-$ dimensional Kenmotsu manifold. On condition that $(\alpha,\beta)=(-1,\beta)$ or $(\alpha,\beta)=(\alpha,0)$, $M$ is viewed as $\xi-$ projectively flat in regard to generalized connection.
\end{theorem}
\begin{corollary}
Let us consider that $M$ is an $n-$ dimensional Kenmotsu manifold, which provides the following expressions: \\
\textbf{(i)} The manifold is suggested to be $\xi-$ projectively flat in regard to $\beta-$ quarter symmetric metric connection.\\
\textbf{(ii)} The manifold is suggested to be $\xi-$ projectively flat in regard to semi-symmetric metric connection.

\end{corollary}
\begin{definition}
Supposing that $M$ is an $n-$ dimensional Kenmotsu manifold, $M$ is suggested to be $\psi-$ projectively flat in regard to the connection  $\overline{D}$ when  $g(\overline{P}(\psi U,\psi V)\phi W,\psi Z)=0$ on $M$.
\end{definition}
By means of (\ref{4.1}), we obtain
\begin{eqnarray}
\overline{P}(\psi U,\psi V)\phi W=\overline{R}(\psi U,\psi V)\phi W-\frac{1}{n-1}\{\overline{S}(\psi V,\psi W)U-\overline{S}(\psi U,\psi W)\psi V\}. \label{fi}
\end{eqnarray}
When (41) and (48) are used, $\psi-$ projectively flatness indicates
\begin{eqnarray}
\overline{K}(\psi U,\psi V,\psi W,\psi Z)&=&\frac{1}{n-1}\{[\overline{S}( V, W)+(n-1)(1+\alpha)]g(\psi U,\psi Z)\ \nonumber\\
&&+[\overline{S}( U, W)+(n-1)(1+\alpha)]g(\psi V,\psi Z)\}. \label{fii}
\end{eqnarray}
When we consider that $\{\nu_{1},\nu_{2},...,\nu_{n-1},\xi\}$ is a local orthogonal basis of the vector fields in $M$ and  that $\{\psi \nu_{1},\psi \nu_{2},...,\psi \nu_{n-1},\xi\}$ is a local orthogonal basis as well, the following is obtained via putting $U=Z=\nu_{i}$ and summing up in regard to $i=1,2,...,n-1$.

\begin{eqnarray}
\sum_{i=1}^{n-1}\overline{K}(\psi \nu_{i},\psi V,\psi W,\psi \nu_{i})&=&\frac{1}{n-1}\{\sum_{i=1}^{n-1}[\overline{S}( V, W)+(n-1)(1+\alpha)]g(\psi e_{i},\psi e_{i})\ \nonumber\\
&&+[\overline{S}(\nu_{i},W)+(n-1)(1+\alpha)]g(\psi V,\psi \nu_{i})\}. \label{6.7}
\end{eqnarray}
We know that
\begin{eqnarray}
\sum_{i=1}^{n-1}g(\psi \nu_{i},\psi \nu_{i})=n-1, \label{6.9}
\end{eqnarray}
and from (\ref{3.19})
\begin{eqnarray}
\sum_{i=1}^{n-1}\overline{S}(\nu_{i},W)g(\psi V,\psi \nu_{i})=\overline{S}( V, W)+(n-1)(1+\alpha)(1-\sum_{i=1}^{n-1}g(\psi V,\psi \nu_{i})). \label{6.10}
\end{eqnarray}
Thus from (\ref{6.9}) and (\ref{6.10}), the equation (\ref{6.7}) becomes
\begin{eqnarray}
\sum_{i=1}^{n-1}\overline{K}(\psi \nu_{i},\psi V,\psi W,\psi \nu_{i})=\frac{n-2}{n-1}\overline{S}( V, W)+(n-1)(1+\alpha). \label{k}
\end{eqnarray}
Moreover we have
\begin{eqnarray}
\sum_{i=1}^{n-1}\overline{K}(\psi \nu_{i},\psi V,\psi W,\psi \nu_{i})=\overline{S}(\psi V,\psi W)-g(\overline{R}(\xi,\psi V)\psi W,\xi) . \label{kk}
\end{eqnarray}
The use of (\ref{3.17}), (\ref{3.19}) and (\ref{kk}), we obtain
\begin{eqnarray}
\sum_{i=1}^{n-1}\overline{K}(\psi \nu_{i},\psi V,\psi W,\psi \nu_{i})=\overline{S}(V,W)+(\alpha+1)\{n-1+g(\psi V,\psi W)+\beta g(\psi V,W)\} \label{kkk}
\end{eqnarray}
Thus, through use of (\ref{k}) and (\ref{kkk}), we have
\begin{eqnarray}
\overline{S}(V,W)=(1-n)(\alpha+1)\{g(V,W)-\eta(V)\eta(W)+\beta g(\psi V,W)\} \label{kkak}
\end{eqnarray}
Accordingly, the following theorem is obtained
\begin{theorem}
The manifold is viewed as generalized $\eta$ Einstein in regard to generalized connection $\overline{D}$ when a Kenmotsu manifold $\psi-$ projectively flat concerning the connection $\overline{D}$.
\end{theorem}
\begin{corollary}
Let us consider that $M$ is an $n-$ dimensional Kenmotsu manifold. When the manifold is $\psi-$ projectively flat In relation to the connection $\overline{\nabla}$, the following expressions are obtained:\\
\textbf{(i)} The manifold is $\eta-$ Einstein manifold in regard to $\beta-$ quarter symmetric connection. \\
\textbf{(ii)} It is Ricci flat regarding generalized symmetric metric connection of type $(-1,\beta)$.\\
\textbf{(iii)} It is Ricci flat concerning semi-symmetric connection.

\end{corollary}

\section{Concircular Curvature Tensor}
When $M$ is considered to be $n-$ dimensional Kenmotsu manifold, the concircular curvature tensor of $M$ in relation the connection $\overline{D}$ can be illustrated through

\begin{eqnarray}
\overline{C}^{*}(U,V)W=\overline{R}(U,V)W-\frac{\overline{r}}{n(n-1)}\{g(V,W)U-g(U,W)V\} \label{coin}
\end{eqnarray}
The use of (\ref{3.2}), (\ref{3.122}) and (\ref{coin}) provides us the following:
\begin{eqnarray}
g(\overline{C}^{*}(U,V)W,\xi)&=&g(C^{*}(U,V)W,\xi)+\{(-3+n)\alpha^{2}+(-4+2n)\alpha\}g(V,W)\eta(U)\ \nonumber\\
&&+\{(3-n)\alpha^{2}+(4-2n)\alpha\}g(U,W)\eta(V)\ \nonumber\\
&&+(\alpha+\alpha^{2})\{g(V,W)\eta(U)-g(U,W)\eta(V)\}\ \nonumber\\
&&+(\beta+\alpha\beta)\{-g(V,\psi W)\eta(U)+g(U,\psi W)\eta(V)\},
 \label{coink}
\end{eqnarray}
where $C^{*}$ is viewed as concircular curvature tensor of $M$ in relation to Levi-Civita connection. When $\overline{C}^{*}=C^{*}$ is considered, we obtain the following by putting $U=\xi$ of the equation (\ref{coink});
\begin{eqnarray}
\{(-2+n)\alpha^{2}+(-3+2n)\alpha\}\{g(V,W)-\eta(V)\eta(W)\}=(\beta+\alpha\beta)g(V,\psi W). \label{etaa}
\end{eqnarray}
The following theorem is obtained through (\ref{Ricci}), (\ref{coink}) and (\ref{etaa})
\begin{theorem}
In a Kenmotsu manifold, the following is attained when concircular curvature tensor is invariant under generalized metric connection
$$\overline{S}(V,W)=S(V,W)-\{(-4+2n)\alpha^{2}+(-6+4n)\alpha\}g(V,W)+\{(2n-4)\alpha^{2}+(3n-5)\alpha\}\eta(V)\eta(W)$$
for every $V,W\in\Gamma(TM)$.
\end{theorem}
\begin{corollary}
Suppose that $M$ is an n- dimensional Kenmotsu manifold. When concircular curvature tensor is invariant under the connection $\overline{\nabla}$, the expressions presented below are obtained:\\
\textbf{(i)} If the manifold is Ricci flat, it is $\eta-$ Einstein regarding generalized symmetric metric connection.\\
\textbf{(ii)} Ricci tensor is invariant  with reference to $\beta-$ quarter symmetric metric connection.

\end{corollary}

\section{Example}

A $3-$ dimensional manifold $M=\{ (x,y,z) \in R^{3}: x\neq 0  \}$ is considered, in which $(x, y, z)$ are regarded as the standard coordinates in $R^{3}$. Suppose that ${\nu_{1},\nu_{2},\nu_{3}}$ are linearly independent global frame on $M$ as presented below
\begin{eqnarray}
\nu_{1}=x\frac{\partial}{\partial z},\;\nu_{2}=x\frac{\partial}{\partial y},\;\nu_{3}=-x\frac{\partial}{\partial x}. \label{5.1}
\end{eqnarray}
Consider that g is a Riemannian metric as presented below
 $$g(\nu_{1},\nu_{2})=g(\nu_{1},\nu_{3})=g(\nu_{2},\nu_{3})=0, g(\nu_{1},\nu_{1})=g(\nu_{2},\nu_{2})=g(\nu_{3},\nu_{3})=1 ,$$
When we consider that the $\eta$ is the $1-$ form represented as $\eta(Z)=g(Z,\nu_{1})$ for every $Z\in TM$ and $\psi$ is the $(1, 1)$ tensor field presented as $\psi \nu_{1}=\nu_{2},\psi \nu_{2}=-\nu_{1}$ and $\psi \nu_{3}=0$, we thereby get $\eta(\nu_{3})=1,\;\psi^{2}Z=-Z+\eta(Z)\nu_{3}$ and $g(\psi Z,\psi W)=g(Z,W)-\eta(Z)\eta(W)$ for all $Z,W\in TM$ through use of linearity of $\psi$ and $g$. Therefore, $\nu_{3}=\xi$ describes an almost contact metric manifold.
Considering that $D$ is the Levi-Civita connection concerning the Riemannian metric $g$, The following is obtained;

\begin{eqnarray}
[\nu_{1},\nu_{2}]=0,\qquad [\nu_{1},\nu_{3}]=\nu_{1},\qquad [\nu_{2},\nu_{3}]=\nu_{2},
\end{eqnarray}
By means of using Koszul formula, the following can be calculated in an easy way
\begin{eqnarray}
D_{\nu_{1}}\nu_{1}=-\nu_{3},\qquad D_{\nu_{1}}\nu_{2}=0. \qquad D_{\nu_{1}}\nu_{3}=\nu_{1}, \nonumber \\
D_{\nu_{2}}\nu_{1}=0, \qquad D_{\nu_{2}}\nu_{2}=-\nu_{3},\qquad D_{\nu_{2}}\nu_{3}=0,\\
D_{\nu_{3}}\nu_{1}=0,\qquad D_{\nu_{3}}\nu_{2}=0,\qquad D_{\nu_{3}}\nu_{3}=0. \nonumber
\end{eqnarray}
The relations presented above remark that
 $(D_{U}\phi)V=g(\phi U,V)\xi-\eta(V)\phi U$, $\nabla_{U}\xi=U-\eta(U)\xi,$ for all $\nu_{3}=\xi $.
 Thus, the manifold $M$ is viewed as a Kenmotsu one with the structure $(\psi,\xi,\eta,g)$. Besides, the non-vanishing components of the tensor can be calculated by using the relations above as presented in the following;
\begin{eqnarray}
R(\nu_{1},\nu_{2})\nu_{1}=\nu_{2},\;R(\nu_{1},\nu_{2})\nu_{2}=-\nu_{1},\,R(\nu_{1},\nu_{3})\nu_{1}=\nu_{3}\nonumber \\
R(\nu_{1},\nu_{3})\nu_{3}=-\nu_{1},\;R(\nu_{2},\nu_{3})\nu_{2}=\nu_{3},\;R(\nu_{2},\nu_{3})\nu_{3}=-\nu_{2} \label{re}
\end{eqnarray}
The equations (\ref{re}) enable us to calculate the components of Ricci tensor in an easy way as in the following: \cite{sukla}

\begin{eqnarray}
S(\nu_{1},\nu_{1})=-2,\;S(\nu_{2},\nu_{2})=-2,\;S(\nu_{3},\nu_{3})=-2
\end{eqnarray}

Now, we can make similar calculations for generalized metric connection. Using (\ref{konnek}) in the above equations, we get
\begin{eqnarray}
\overline{D}_{\nu_{1}}\nu_{1}=-(1+\alpha)\nu_{3},&  \overline{D}_{\nu_{1}}\nu_{2}=0. &  \overline{D}_{\nu_{1}}\nu_{3}=(1+\alpha)\nu_{1}, \nonumber \\
\overline{D}_{\nu_{2}}\nu_{1}=0, & \overline{D}_{\nu_{2}}\nu_{2}=-(1+\alpha)\nu_{3},& \overline{D}_{\nu_{2}}\nu_{3}=\alpha \nu_{2},\label{na} \\
\overline{D}_{\nu_{3}}\nu_{1}=-\beta \nu_{2},& \overline{D}_{\nu_{3}}\nu_{2}=\beta \nu_{1}, &\overline{D}_{\nu_{3}}\nu_{3}=0. \nonumber
\end{eqnarray}
By means of (\ref{na}), we can make calculations of the components of curvature tensor concerning generalized metric connection as follows:

\begin{eqnarray}
\overline{R}(\nu_{1},\nu_{2})\nu_{1}=(1+\alpha)^{2}\nu_{2},\;& \overline{R}(\nu_{1},\nu_{2})\nu_{2}=-(1+\alpha)^{2}\nu_{1},\nonumber \\
\overline{R}(\nu_{1},\nu_{3})\nu_{1}=(1+\alpha)\nu_{3}      & \overline{R}(\nu_{1},\nu_{3})\nu_{3}=(1+\alpha)(\beta \nu_{2}-\nu_{1}),\;\nonumber  \\
\overline{R}(\nu_{2},\nu_{3})\nu_{2}=(1+\alpha)\nu_{3},\;& \overline{R}(\nu_{2},\nu_{3})\nu_{3}=-(1+\alpha)(-\beta \nu_{1}+ \nu_{2})\ \label{r} \\
\overline{R}(\nu_{3},\nu_{2})\nu_{1}=-(1+\alpha)\beta \nu_{3},\;& \overline{R}(\nu_{3},\nu_{1})\nu_{2}=(1+\alpha)\beta \nu_{3},\, \nonumber
. &
\end{eqnarray}
Through (\ref{r}), The components of Ricci tensor mentioned above are presented below:

\begin{eqnarray}
\overline{S}(\nu_{1},\nu_{1})=-(1+\alpha)(2+\alpha),\;& \overline{S}(\nu_{2},\nu_{2})=-(1+\alpha)(2+\alpha),\ \nonumber \\
\overline{S}(\nu_{3},\nu_{3})=-2(1+\alpha),\;& \overline{S}(\nu_{2},\nu_{1})=-(1+\alpha)\beta. \ \label{s}
\end{eqnarray}
(\ref{r}) and (\ref{s}) are verified through the equations (\ref{eg}) and (\ref{Ricci}), respectively. Moreover, the scalar curvature regarding the Levi-Civita connection and generalized-metric connection are $r=-6$ and $\overline{r}=-2(1+\alpha)(3+\alpha)$.
The following expressions are obtained through (\ref{4.1}):
\begin{eqnarray}
\overline{P}(\nu_{1},\nu_{3})\nu_{3}=(1+\alpha)\beta \nu_{2},\;\overline{P}(\nu_{2},\nu_{3})\nu_{3}=-(1+\alpha)\beta E_{1},\;\overline{P}(\nu_{1},\nu_{2})\nu_{3}=0 .\label{pro}
\end{eqnarray}
which verify Theorem $5$, Theorem $6$ and Theorem $11$.

\noindent Mohammed Ali Alghamdi

\noindent Department of Mathematics, Faculty of Science,  King Abdulaziz University

\noindent Jeddah 21589, Saudi Arabia

\noindent E-mail: proff-malghamdi@hotmail.com

\noindent O\u{g}uzhan Bahad\i r

\noindent Department of Mathematics, Faculty of Arts and Sciences, \ K.S.U

\noindent Kahramanmaras, TURKEY

\noindent Email: oguzbaha@gmail.com

\end{document}